\documentclass[12pt]{article}
\usepackage{amssymb}

\newcommand{\A}{{\bf \cal A}}

\newcommand{\E}{{\Bbb E}}
\newcommand{\F}{{\cal F}}

\newcommand{\R}{{\Bbb R}}
\newcommand{\half}{{  {1\over 2}  }}

\newtheorem{theorem}{Theorem}[section]
\newtheorem{proposition}[theorem]{Proposition}
\newtheorem{lemma}[theorem]{Lemma}
\newtheorem{corollary}[theorem]{Corollary}

\def\exp{{\rm e}}

\begin{document}

\title{A class of integration by parts formulae in stochastic analysis I}
\author{K. D.   Elworthy \ and Xue-Mei Li \\Mathematics Institute\\
University of Warwick\\Coventry CV4 7AL,U.K.}

\date{}
\maketitle

\section{Introduction}


Consider a Stratonovich stochastic differential equation
\begin{equation}\label{1}
dx_t=X(x_t)\circ dB_t +A(x_t)dt
\end{equation}
with $C^\infty$ coefficients on a compact Riemannian manifold $M$, with 
associated differential generator $\A=\half \Delta_M+ Z$ and solution flow
$\{\xi_t: t\ge 0\}$ of random smooth diffeomorphisms of $M$. Let
$T\xi_t:
TM\to TM$ be the induced map on the tangent bundle of $M$ obtained by
differentiating $\xi_t$ with respect to the initial point.  Following  an
observation by  A. Thalmaier we  extend the basic formula of \cite{EL-LI}
to obtain
\begin{equation}\label{2}
\E dF\left(T\xi_\cdot\left(h_\cdot\right)\right)=
\E F\left(\xi_\cdot(x)\right) \int_0^T 
\left<T\xi_s \left (\dot h_s \right), X\left(\xi_s(x)\right)dB_s \right>
\end{equation}
where $F\in \F C_b^\infty(C_x(M))$, the space of smooth cylindrical functions
on the space $C_x(M)$ of  continuous paths
 $\gamma:[0,T]\to M$ with $\gamma(0)=x$, $dF$ is its derivative,
  and $h_\cdot$ is a suitable
 adapted process with sample paths in the Cameron-Martin space
 $L_0^{2,1}([0,T]; T_xM)$. Set $\F_t^x=\sigma\{\xi_s(x): 0\le s\le t\}$.
Taking conditional expectation with respect to $\F_T^x$, formula
(\ref{2}) yields integration by parts formulae on $C_x(M)$ of the form
\begin{equation}\label{3}
\E dF(\gamma)(\bar V^h)=\E F(\gamma)\delta\overline{ V}^h(\gamma)
\end{equation}
where $\bar V^h$ is the vector field on $C_x(M)$
$$ \bar V^h(\gamma)_t =\E\left\{T\xi_t(h_t) \left| \xi_\cdot(x)=
\gamma\right. \right\}$$
and 
$\delta\overline{V}^h: C_x(M)\to \R$ is given by
$$\delta\overline{V}^h(\gamma)=\E\left\{\int_0^T<T\xi_s(\dot h_s),
 X(\xi_s(x))dB_s> \left|  \xi_\cdot(x)=\gamma\right. \right\}. $$

\bigskip

When $h_\cdot$ is adapted to $\F_\cdot^x$ results from \cite{EL-LJ-LI1}
extending \cite{EL-YOR} give explicit expressions for $\bar V^h$ and
 $\delta \bar V^h$ in terms of the Ricci curvature of the LeJan-Watanabe
connection associated to (\ref{1}). Equation (\ref{3}) then reduces
to a Driver's integration by parts formula, Theorem \ref{3.3} below,
but no hypothesis of torsion skew symmetry of the connection is required:
 the integration by parts formulae follow for the adjoint of any
metric connection.
  In particular for any such connection there is a Hilbert ``tangent
space'' of ``good'' directions obtained by parallel translation of 
the Cameron-Martin space of paths in $T_xM$. (In fact
it is the ``Ricci flow'' or ``Dohrn-Guerra parallel translation''
(see Nelson \cite{Nelson84}), leading to the ``damped
gradient'' (\cite{Fang-Malliavin93}) which occurs more naturally.) 
However, in Remark 2.4, we show that in this case
$\bar V_h$ is in the class for which integration by parts
 formulae are known, so that the results of 2.3, 3.3, 3.5 are not
claimed to be new in substance.

\bigskip

Although this filtering out of the extraneous noise gives intrinsic
results comparable to those of Driver \cite{Driver92},
 this viewpoint throws away a lot of the structure we have.
Moreover integration by parts formulae such as (\ref{2}) should
have some connection with quasi-invariance properties of flows
 associated to the vector fields. Flows for the $\bar V^h$ on $C_x(M)$
 do not appear to be  easy to analyse in general.
However in \S 3 we  show that in the context of Diff M valued processes
 there are very natural flows associated and (\ref{2}) has a rather
 natural geometric interpretation. This leads to another
elementary proof of (\ref{2}) and in Theorem \ref{3.4} we use this
method to obtain integration by parts formulae for the free path space.

\bigskip

There are at least 3 proofs of (\ref{2}). The first given here
is via It\^o's formula and elementary martingale calculus (it
requires $F$ to be cylindrical), the second given here is based
 on the Girsanov-Maruyama theorem (and works for more general
$F$), and a third method would be to deduce it from the standard
integration by parts formula on Wiener space applied to the
 functional $F\circ \xi$, c.f. \cite{Bismut81b}.
 Indeed this work was stimulated by
D. Bell and  D. Nualart pointing out that this third approach could be 
used to deduce the basic formula of \cite{EL-LI}. The point
made (and carried out) in \cite {ELflow} and \cite{EL-LI}
that the first approach  can be applied directly to 'Ricci flows' 
instead of derivative flows to give intrinsic formulae without
stochastic flows, 
also needs to be emphasized: see also \cite{Stroock-Zeitouni}.
As such it gives the details of how `Bismut's formula' (essentially
integration by parts when $F$ is a function of paths evaluated
at just one time $t$) leads to the full integration by parts formula.

There are also now many proofs of Driver's results for $C_x(M)$ and 
for the free path space and their extensions. See \cite{Hsu95},
\cite{Enchev-Stroock95},  \cite{Leandre-Norris} (with a very 
concise proof), \cite{Airault-Malliavin},
 \cite{Aida95}, and  \cite{Cruzeiro-Malliavin}. 

\bigskip

\noindent{\bf Acknowledgment:}
This research was supported by SERC grant GR/H67263 and stimulated and
helped by our contacts with A. Thalmaier.

\section{The integration by parts formula from finite dimensional
manifolds to path spaces}

In this section we deduce by induction an integration by parts formula on
the path space from a  formula on the base  manifold $M$.
The key is to obtain formula (\ref{formula2}) for $M$.

\bigskip

Let $h: \Omega\times [0,T] \to T_{x}M$ be an adapted process with 
$h(\omega): [0,T]\to T_{x}M$ in $L^{2,1}$ for almost all $\omega$.

\begin{lemma}\label{2.1}
If $h: \Omega\times [0,T]\to T_{x}M$ is adapted,  $L^{2,1}$ for
 a.s. $\omega$ and
$\left(\int_0^T|\dot h_s|^2 ds\right)^{1/2}\in L^{1+\epsilon}$ for
some $\epsilon>0$.  Then for $t< T$,

\begin{equation} 
\label{new4}
\begin{array}{l}
\E\left\{\int_0^t <T\xi_s(\dot h_s), X(\xi_s(x))dB_s> \left |
\xi_T(x)\right.  \right\}\\
=\E\left\{ \int_t^T <T\xi_s(-), X(\xi_s(x))dB_s> {h_t-h_0\over T-t}
\left| \xi_T(x)\right.\right\}.
\end{array}\end{equation}

If furthermore $h_\cdot$ is non-random then for $t\le T$,
\begin{equation}\label{4}\begin{array}{l}
\E\left\{\int_0^t  <T\xi_s(\dot h_s), X(\xi_s(x))dB_s>
\left| \xi_T(x)\right. \right\}\\
=\E\left\{ \int_0^t<T\xi_s(-), X(\xi_s(x))dB_s>\left({h_t-h_0\over t}\right)
\left| \xi_T(x)\right. \right\}.
\end{array}
\end{equation}
\end{lemma}

\noindent
{\bf Proof.}
First by the Burkholder-Davis-Gundy inequality, for some constant $c_1$,
\begin{eqnarray*}
&&\E\left|\int_0^T <T\xi_s(\dot h_s), X(\xi_s(x))dB_s>\right|
 \le c_1 \E \left(\int_0^T |T\xi_s(\dot h_s)|^2 ds\right)^{1\over 2}\\
&&\le c_1\left(\E \sup_{0\le s\le T} |T_x\xi_s|^{1+\epsilon \over
  \epsilon}\right) ^{\epsilon\over 1+ \epsilon}
\left[\E\left(\int_0^T |\dot h_s|^2 ds\right)^{1+\epsilon\over
  2}\right]^{1\over 1+\epsilon}.
\end{eqnarray*}
This is finite since $\sup_{0\le s\le t} |T_x\xi_s|\in L^q$ for all
$1\le q<\infty$, e.g. see \cite{application}.  Moreover, since
the adapted processes in
 $L^\infty(\Omega, \F, {\Bbb P}; C^1\left([0,T];T_xM\right))$
are dense in the subspace of adapted processes in 
$L^{1+\epsilon}\left(\Omega, \F, {\Bbb P};
L^{2,1}\left([0,T];T_xM\right) \right)$,   this estimate allows us to
 assume that $h$ belongs to the former space.

\bigskip

Set 
 $M_t=\int_0^t<T_x\xi_s(-), X(\xi_s(x))dB_s>$.
Then $\{M_\cdot\}$ is a $T_{x}^* M$ valued local martingale. If
$0=t_0<t_1<\dots<t_l=t$ is a partition of $[0,t]$,
 $\Delta_jt=t_{j+1}-t_j$, and  $\Delta_jM=M_{t_{j+1}}-M_{t_j}$, then
\begin{equation}\label{5}
\sum_{j=1}^{l-1}\Delta_jM( \dot h_{t_j}) \to
\int_0^t \dot h_sdM_s
=\int_0^t<T\xi_s(\dot h_s),  X(\xi_s(x))dB_s>
\end{equation}
and the convergence is in $L^1$.

 On the other hand if $v_0\in T_xM$
and $P_t$ is the probabilistic semigroup associated to the
S.D.E. and $f$ a bounded measurable function then
\begin{equation}\label{formula}
d(P_Tf)(v_0)={1\over T}\E f(\xi_T(x)) 
\int_0^T\left<T\xi_s(v_0), X(\xi_s(x))dB_s\right>.
\end{equation}
See \cite{EL-LI}. However by an observation of Thalmaier:
the same proof shows that  for
any $r, h\in [0,T]$ with $h>0$ and $r+h\le T$
$$d(P_Tf)(v_0)={1\over h}\E f(\xi_T(x)) \int_r^{r+h}
\left<T\xi_s(v_0), X(\xi_s(x))dB_s\right>$$
c.f. \cite{Stroock-Zeitouni}. From these two formulae
we obtain:
\begin{equation}
\label{7}
\begin{array}{l}
\E\left\{{1\over T}\int_0^T <T\xi_s(v_0), X(\xi_s(x))dB_s>
\left| \xi_T(x) \right.\right\}\\
=\E\left\{{1\over h}\int_r^{r+h} <T\xi_s(v_0), X(\xi_s(x))dB_s>
\left| \xi_T(x)\right.\right\}.
\end{array}
\end{equation}
For any $0\le r\le T$, let $\{\xi_s^r(x): r\le s\le T, x\in M\}$
be the solution flow to (\ref{1}) starting from $x$ at time $r$.
The flow $\xi_\cdot^r$ can be taken to be adapted to a filtration
$\{\F_s^r: r\le s\le T\}$ independent of $\F_r$, and then
we have $\xi_s^r\xi_r=\xi_s$, almost surely, $r\le s\le T$.
>From this, time homogeneity, and (\ref{7}),

\begin{eqnarray*}
&&\E\left\{\sum_{j=1}^{l-1}
 \Delta_jM(\dot h_{t_j})\left|\xi_T(x)\right.\right\}\\
&&=\E   \left\{ \sum_{j=1}^{l-1}  \Delta_jt {1\over \Delta_j t}
\int_{t_j}^{t_{j+1}} 
\left<    T\xi_s^{t_j}\left(T\xi_{t_j}  \left(\dot h_{t_j})\right)\right),
X\left(\xi_s^{t_j}\left(\xi_{t_j}(x)\right)\right)dB_s\right>  \,  \left|
\xi_T^{t_j}(\xi_{t_j}(x))    \right.\right\}\\
&&=\E\left\{\sum_{j=1}^{l-1}  \Delta_j t {1\over T-t} \int_t^T
\left<T\xi_s^{t_j}\left(T\xi_{t_j}   (\dot
h_{t_j})\right),
X\left(\xi_s^{t_j}\left(\xi_{t_j}(x)\right)\right)dB_s\right>\, \left|
\xi_T^{t_j}(\xi_{t_j}(x))\right.\right\}\\
&&=\E\left\{\sum_{j=1}^{l-1} \Delta_{j}t {1\over T-t}
\int_t^T<T\xi_s(\dot h_{t_j}), X(\xi_s(x))dB_s>
\left|\xi_T(x)\right.\right\}\\
&&\to \E\left\{ \int_{t}^{T}<T\xi_s(-), X(\xi_s(x))dB_s> {h_t-h_0\over T-t}
\left| \xi_T(x)\right.\right\}.\\
\end{eqnarray*}

Comparing with (\ref{5}) this gives the first required identity. When
$h_\cdot$ is non-random the second follows immediately from
(\ref{7}).
\hfill \rule{2mm}{3mm}

\bigskip

\noindent
{\bf Remark:}

As in \cite{Stroock-Zeitouni} a further modification is possible
replacing (\ref{7}) by:

\begin{eqnarray*}
&&{1\over T}\E\left\{\int_0^T<T\xi_s(v_0), X(\xi_s(x))dB_s>\left|
 \xi_T(x)\right.\right\}\\
&&= {1\over \int_0^T \Psi(r)dr}
\E\left\{\int_0^T \Psi(s) <T\xi_s(v_0), X(\xi_s(x))dB_s>\left|
 \xi_T(x)\right.\right\}
\end{eqnarray*}
for $\Psi: [0,T]\to \R$ integrable  with $\int_0^T \Psi(r)dr\not =0$.
The argument leads to, for non-random $h$,
\begin{equation}\label{5000}
\begin{array}{l}
\E\left\{\int_0^t  <T\xi_s(\dot h_s), X(\xi_s(x))dB_s>
\left| \xi_T(x)\right. \right\}\\
=\E\left\{\int_0^t  \Psi(s)<T\xi_s(-), X(\xi_s(x))dB_s>
\left(  {h_t-h_0 \over \int_0^t\Psi(r) dr}\right)
\left| \xi_T(x)\right. \right\}.
\end{array}
\end{equation}

\bigskip

\begin{corollary}
 Under the conditions of the lemma, for any  $C^1$ function
$f: M\to \R$,
\begin{equation}
\label{formula2}
\E f\left(\xi_T(x)\right)\int_0^T<T\xi_s(\dot h_s), X(\xi_s(x))dB_s>
=\E df\left(T\xi_T\left(h_T-h_0\right)\right).
\end{equation}
\end{corollary}

\noindent{\bf Proof.}
First by the composition property of solution flows,
\begin{eqnarray*}
&& \E\left\{
\int_{t}^{T}<T\xi_s(-), X(\xi_s(x))dB_s> {h_t-h_0\over T-t}
\left| \xi_T(x)\right.\right\}\\
&&=\E\left\{\int_{t}^{T}<T\xi_s^t(-),
X(\xi_s^t\left(\xi_t(x)\right))dB_s>{ T\xi_t(h_t-h_0)\over T-t}
\left| \xi_T^t\left(\xi_t(x)\right) \right.\right\}. \\
\end{eqnarray*}

As in the proof of the lemma, (\ref{new4}) yields
\begin{eqnarray*}
&&\E f(\xi_T(x))\int_0^t <T\xi_s(\dot h_s), X(\xi_s(x)) dB_s>\\
&&=\E f(\xi_T^t(\xi_t(x))\int_t^T
 \left<T\xi_s^t(-), X(\xi_s^t(\xi_t(x))dB_s\right>\,
{T\xi_t(h_t-h_0)\over T-t}\\
&&=\E\left\{dP_{T-t}(f) \left(T\xi_t(h_t-h_0)\right)\right\}
\end{eqnarray*}
by \cite{EL-LI}, since $\F_\cdot^t$ is independent of $\F_t$.
 Now let $t$ increase to $T$ and  the required result follows.
\hfill \rule{2mm}{3mm}

\bigskip

Next consider a cylindrical function $F$ on $C_{x}(M)$, the space of
continuous paths with base point $x$. Write
$$F(\gamma_\cdot)=f(\gamma_{t_1}, \dots, \gamma_{t_k}),$$
for $(t_1, \dots, t_k)\in [0,T]^k$, $\gamma\in C_x(M)$ and 
 $f$ a smooth function on $M^k$. Suppose $h_0=0$ and consider 
the tangent vector
field $V^h(\xi_\cdot(x))$ along $\{\xi_t(x): 0\le t\le T\}$ 
on $C_{x}(M)$ given by
$$V^h(\xi_\cdot)_t=T_{x}\xi_t(h_t).$$
Then 
\begin{equation}
dF(V^h(\xi_\cdot))=\sum_{j=1}^k d^j f_{\xi_{\underline t}}
\left(V^h(\xi_\cdot)_{t_j}\right).
\end{equation}
Here $\xi_{\underline t}=(\xi_{t_1},\dots, \xi_{t_k})$ and
 $d^jf$ is the partial derivative of $f$ in the $j$th direction.

Let 
$$\delta V^h(\xi_\cdot)
=\int_0^{T}<T_{x}\xi_s(\dot h_s), X(\xi_s(x))dB_s>.$$

\begin{theorem} Let $h: [0,T]\times \Omega\to T_{x}M$ be an adapted
 stochastic process with  almost surely all $h(\omega)\in L_0^{2,1}$
and $\E\left(\int_0^T |\dot h_s|^2ds\right)^{1+\epsilon\over
  2}<\infty$ for some $\epsilon>0$.  Then
\begin{equation}\label{8}
\E dF(V^h(\xi_{\cdot}))=\E F(\xi_{\cdot}(x))\delta V^h(\xi_{\cdot}).
\end{equation}
\end{theorem}

\noindent{\bf Proof.} We prove by induction on $k$. When $k=1$, this
is just  (\ref{formula2}), the formula for functions. 
Let $\Omega=C_0([0, T]; \R^n)$ be  the canonical probability space.
We set $\Omega_1=C_0([0, t_1]; \R^n)$ and
 $\Omega_2=C_0([t_1, T]; \R^n)$. There  is then the standard 
decomposition of filtered spaces
\begin{eqnarray*}
&&\{\Omega, \F, \F_t, 0\le t\le T, {\Bbb P}\}\\
=&&\{\Omega_1, \F, \F_t, 0\le t\le t_1, {\Bbb P_1}\} 
\times \{\Omega_2, \F, \F_t^{t_1}, t_1\le t\le T, {\Bbb P_2}\}
\end{eqnarray*}
in the sense that $\F_t=\F_t*\Omega_2$ if $t\le t_1$,
and $\F_t=\F_{t_1}*\F_t^{t_1}$ if $t\ge t_1$. As before let
 $\xi_t^{t_1}(y_0), t_1\le t\le T, y_0\in M$ be the solution flow 
to (\ref{1}) starting at time $t_1$, i.e.  $\xi_{t_1}^{t_1}(y_0)=y_0$.
We will consider it as a function of $\omega_2 \in \Omega_2$, 
adapted to $\F^{t_1}_\cdot$, while $\{\xi_t: 0\le t\le t_1\}$
will be considered on $\Omega_1$, and  $\{\xi_t: t_1\le t\le T\}$
on $\Omega_1\times \Omega_2=\Omega$. The composition property for flows
gives
$$\xi_t^{t_1}\left(\xi_{t_1}(x, \omega_1), \omega_2\right)
=\xi_t\left(x, (\omega_1, \omega_2)\right), \hskip 10 pt \hbox{each }
t_1\le t\le T, a.s.$$

Assume the required result holds for cylindrical functions depending
on $k-1$ times, some $k\in \{2, 3\dots\}$.  Take  $y_0\in M$ and
define  $f_1^{y_0}: M^{k-1}\to \R$ and $F_1^{y_0}: \Omega_2\to \R$  by:
$$ f_1^{y_0}(x_1, \dots, x_{k-1})=f(y_0, x_1, \dots,  x_{k-1})$$
and
$$ F_1^{y_0} (\omega_2)=f(y_0, \xi_{t_2}^{t_1}(y_0, \omega_2), \dots,
 \xi_{t_k}^{t_1}(y_0, \omega_2)).$$

Take $h_\cdot^1: \Omega_2\to L_0^{2,1}\left([t_1, T]; T_{y_0}M\right)$,
 adapted to $\F^{t_1}_\cdot$, and with 
$\E\left(\int_{t_1}^T |\dot h_s^1|^2 ds\right)^{1+\epsilon\over 2}$
 finite. By time homogeneity our  inductive hypothesis gives

\begin{equation}\label{1000}
\begin{array}{l}
\sum_{j=2}^k\int_{\Omega_2}
 d^j f\left(y_0, \xi_{t_2}^{t_1}(y_0, \omega_2), \dots,
 \xi_{t_k}^{t_1}(y_0, \omega_2)\right)  
\left(T\xi_{t_j}^{t_1}(h_{t_j}^1(\omega_2), \omega_2)\right)
 d{\Bbb P_2}(\omega_2)\\
=\int_{\Omega_2}
f\left(y_0, \xi_{t_2}^{t_1}(y_0, \omega_2), \dots,
 \xi_{t_k}^{t_1}(y_0, \omega_2)\right) \times\\
 \hskip14pt \int_{t_1}^{T} \left<T\xi_{r}^{t_1}(\dot h^1_r(\omega_2),\omega_2),
 X(\xi^{t_1}_r(y_0,\omega_2))d B_r(\omega_2)  \right> 
 d{\Bbb P_2}(\omega_2) .
\end{array}\end{equation}

Now for $\omega_1\in \Omega_1$ (outside of a certain measure zero set)
 we can take  $y_0=\xi_{t_1}(x_0, \omega_1)$ and 
 $$h_t^1(\omega_2)
=T\xi_{t_1}\left(h_{t}\left(\omega_1, \omega_2\right)-h_{t_1}(\omega_1),
 \omega_1\right).$$ 
Then, for almost all $\omega_1\in \Omega_1$, we have $h_\cdot^1$
adapted to $\F_\cdot^{t_1}$.
Substitute this in (\ref{1000}).  Using the composition property,
 and then integrating over $\Omega_1$ yields

\begin{equation}\label{1001}
\begin{array}{l}
\sum_{j=2}^k \E d^j f(\xi_{\underline t}) 
\left( T\xi_{t_j}(h_{t_j}-h_{t_1}) \right)\\
= \E f(\xi_{\underline t}(x)) \int_{t_1}^T
\left< T\xi_r(\dot h_r), X(\xi_r(x))dB_r\right>.
\end{array}\end{equation}


On the other hand we can define  $g: M\to \R^1$ by

$$g(x)=\int_{\Omega_2} f\left(x, \xi_{t_2}^{t_1}(x, \omega_2),
 \dots, \xi_{t_k}^{t_1}(x, \omega_2)\right)$$
and apply  formula (\ref{formula2}) to $g$  to obtain:

$$\int_{\Omega_1} dg(T\xi_{t_1}(h_{t_1}))d{\Bbb P_1}(\omega_1)
=\int_{\Omega_1}  g(\xi_{t_1}(x))\int_0^{t_1}
\left <T\xi_r(\dot h_{r})), X(\xi_r(x_0))dB_r\right>d{\Bbb P_1}(\omega_1).$$
But note that 
$$\int_{\Omega_1} d g(T\xi_{t_1}(h_{t_1}))d{\Bbb P_1}(\omega_1)
=\sum_{j=1}^k \E d^kf_{\xi_{\underline t}}(T \xi_{t_j}(h_{t_1})),$$
and therefore

\begin{equation}\label{9}
\sum_{j=1}^k \E d^jf_{\xi_{\underline t}}(T\xi_{t_j}(h_{t_1}))
=\E f(\xi_{\underline t}) \int_0^{t_1}
 \left< T\xi_r(\dot h_r), X(\xi_r(x))dB_r \right> 
\end{equation}

Adding  (\ref{1001}) we arrive at (\ref{8}):
\begin{eqnarray*}
\sum_{j=1}^k \E d^jf_{\xi_{\underline t}}(T\xi_{t_j}(h_{t_j}))
&=&\E f(\xi_{\underline t}(x))
 \int_0^{T} \left< T\xi_r(\dot h_r), X(\xi_r(x))dB_r \right>. \\
\end{eqnarray*}
\hfill\rule{2mm}{3mm}

{\bf B.} Let $\tilde \nabla$ be a metric connection for the manifold
$M$ with torsion $T$, and $\tilde \nabla^\prime$ its adjoint connection
 defined by
$$\tilde \nabla^\prime_{V_1}V_2=\tilde \nabla_{V_1}V_2-T(V_1,V_2).$$
Here $V_1, V_2$ are  vector fields. Let $\tilde R$ be the curvature
tensor
of $\tilde \nabla$ and define $\tilde{Ric}^\#:
TM\to TM$ by $\tilde{Ric}^\#(v)=\hbox{trace } \tilde{R}(v,-)-$.
 If $\{x_s\}$ is  a diffusion
on $M$ with generator 
$\half {\rm trace} \tilde \nabla{\rm grad}+L_Z$  denote by $\tilde //_s$
 the parallel transport along $\{x_s\}$, and $\{\tilde B_s: 0\le s\le t\}$
the martingale part of the anti-development of $\{x_s: 0\le s\le t \}$
 using $\tilde //_s$,
 a Brownian motion on $T_{x_0}M$. Let $v_s=\tilde W_s^Z(v_0)$ be the
 solution to 
$${\tilde D^\prime \over \partial s}v_s=
-\half \tilde{\rm  Ric}^\#(v_s)+\tilde\nabla Z(v_s)$$
starting from $v_0\in T_{x_0}M$. Here $\tilde D^\prime $ denotes the
 covariant differentiation along the paths of $\{x_t\}$ using the
 adjoint connection.
 We will show that 
 (\ref{8}) implies  Driver's integration by parts formula.
However we do not need to assume $\tilde \nabla^\prime$ (or
equivalently $\tilde \nabla$) is torsion skew symmetric.

\begin{corollary}\label{2.3}
 Let $F$ be a cylindrical function on $C_{x_0}(M)$.  
 Suppose $h: [0,T]\times \Omega\to T_{x_0}M$ is adapted to
 the filtration  of $\{x_s: 0\le s<\infty\}$ 
and such that  $h(\omega)$ is in $ L_0^{2,1}$
 for  almost all $\omega$  and $h\in L^{1+\epsilon}\left(\Omega, \F,
 {\Bbb P}; L_0^{2,1}([0,T]; T_{x_0}M)\right)$ for some $\epsilon>0$. Then
\begin{equation}\label{1500}
\E dF(\tilde W_\cdot^Z(h_\cdot))= \E F(\xi_{\cdot}(x_0))
\int_0^T<\tilde W_s^Z(\dot h_s), \tilde{//}_sd\tilde B_s>.
\end{equation}
When $\tilde \nabla^\prime$ is metric for some Riemannian metric
on $M$, it suffices to have $h\in L^1\left(\Omega, \F, {\Bbb P};
L_0^{2,1}([0,T])\right)$.
\end{corollary}

\noindent{\bf Proof.}
By a result of \cite{EL-LJ-LI1} we can choose $X$ such that
$\tilde \nabla$ equals the Le Jan-Watanabe connection induced from
the stochastic differential equation
$$dx_t=X(x_t)\circ dB_t+Z(x_t)dt$$
and the solution flow $\{\xi_\cdot(x)\}$ has generator
$\half {\rm trace} \tilde \nabla{\rm grad}+L_Z$ (c.f. Corollary 3.4
of \cite{EL-LJ-LI1}). Moreover 
 the conditioned process of the derivative flow $T\xi_t(v_0)$ with
 respect to the natural filtration of $\{\xi_\cdot(x_0)\}$ is given by 
 $\{\tilde W_\cdot^Z(v_0)\}$:
$$\E \{ T\xi_t(v_0)\left| \right. \F_T^{x_0}\}=\tilde W_t^Z(v_0),$$
by Theorem 3.2 of \cite{EL-LJ-LI1}  extending \cite{EL-YOR}.
The result follows since $\tilde B_t$
equals $\int_0^t \tilde{//}_s^{-1}X(\xi_s(x_0))dB_s$.

If $\tilde\nabla^\prime$ is metric for some Riemannian metric
then $\sup_{0\le s\le t}|\tilde W_s^Z|$ is in 
$L^\infty\left(\Omega, \F, {\Bbb P}\right)$ and so the
Burkholder-Davis-Gundy
inequality used as in the proof of Lemma 2.1 allows us to take
$\epsilon=0$.
\hfill \rule{2mm}{3mm}

\bigskip

\noindent{\bf Remarks 2.5.} 
 (i). Let $S: TM\times TM\to TM$ be a tensor fields of type (1,2), and
let $\nabla$ refer to the Levi-Civita connection of $M$. Then,
by \cite{Kobayashi-NomizuII} p.146, a connection $\tilde \nabla$
can be defined by 
$$\tilde \nabla_{V_1}(V_2)=\nabla_{V_1}(V_2)+S(V_1,V_2)$$
for vector fields $V_1$, $V_2$. and all linear connections on
$M$ can be obtained this way. It is easy to see that $\tilde \nabla$
is metric if and only if
$$<S(W,U), V>=-<U, S(W,V)>$$
for all vector fields $U$, $V$, $W$, i.e. if and only
if $S(W,-)$ is skew symmetric. On the other hand the adjoint
connection is given by 
$$\tilde\nabla^\prime_{V_1}(V_2)=\nabla_{V_1}(V_2)+S(V_2, V_1)$$
so that it is torsion skew symmetric if also $S(-,W)$ is skew
symmetric. In terms of the Levi-Civita connection
our vector fields $\bar V^h$  for which the integration by parts
formula hold therefore satisfy an equation of the form

$$D \bar V_t^h= -S(\bar V_t^h, \circ dx_t)+
\Lambda_t(\bar V_t^h)dt+W_t^h(\dot h_t)dt+ \nabla A(\bar v_t^h)dt$$
where $\Lambda_t$ is linear (also depending on $S$). In particular
they are ``tangent processes'' in the sense proposed by Driver,
for which integration by parts formulae are known: see
\cite{Driver95},  \cite{Cruzeiro-Malliavin}, \cite{Airault-Malliavin}, and
\cite{Aida95}, \cite{Driver95II}.

(ii) For cylinder functions depending on one time only such
integration by parts formulae go back to Bismut \cite{Bismut84}.

\section{Geometric intepretation and a shorter proof}

{\bf A.}  
    The processes $T_x\xi_t(h_t)$
 cannot strictly speaking be considered as tangent vectors or 
vector fields on $C_x(M)$. In some sense they form tangent vectors at
 $\xi_\cdot(x,-)$ to the space of processes (or semi-martingales)
$$[0,T]\times \Omega \to M$$
since $T_x\xi_t (h_t(\omega), \omega)\in T_{\xi_t(x,\omega)}M$ for
$(t,\omega)\in[0,T]\times \Omega$ or equivalently as 'tangent vectors'
to the space of random variables 
$$\Omega \to C_x(M)$$
at $\omega \mapsto \xi_\cdot(x,\omega)$. However c.f. \cite{Driver92}
there is still
 no natural associated flow. In fact the most natural interpretation
 takes into account the variable $x$ and replaces $C_x(M)$ by
 $P_{id}$Diff$M$ the space of paths on the diffeomorphism
group of $M$, as we now describe.

\bigskip

 Let Diff$M$ be the space of $C^\infty$ diffeomorphisms of $M$.
We can consider it with a rather formal differential structure or 
if the reader prefers it can be replaced by a suitable Sobolev
space of diffeomorphisms, to give a Hilbert manifold (as in 
\cite{ELbook} following \cite{Ebin-Marsden70}). In any case the tangent space
 $T_\alpha(\hbox{Diff} M)$ will be identified with all vector fields
on $M$ over $\alpha$ i.e. smooth $v: M\to TM$ such that
 $v(x)\in T_{\alpha(x)}M$ for all $x\in M$. If $P$Diff$M$ refers
to continuous paths $\phi:[0,T]\to \hbox{Diff}M$ with $\phi(0)=id_{M}$ 
then $T_\phi P \hbox{Diff} M$ will be identified with continuous 
$v:[0,T]\to T\hbox{Diff} M$ vanishing at $t=0$, 
 such that $v(t)\in T_{\phi(t)} \hbox{Diff} M$,
or equivalently $v: [0,T]\times M\to TM$ with $v(t)(x)\in T_{\phi(t)(x)}M$.

\bigskip

{\bf B.}
 Given our S.D.E. (\ref{1}) now take $h\in L_0^{2,1}([0,T];\R^n)$.
There is $X^{h_\cdot}$, the time dependent vector field $X(\cdot)(h_t)$
on $M$.  From this we obtain a field $U^h$ on $P \hbox{Diff}M$ by

\begin{equation}
U^h(\phi)_t(x)=T_x\phi_t(X(x)h_t).
\end{equation}
This is just the left invariant vector field on $P\hbox{Diff}M$
 corresponding to $X^{h_\cdot}\in T_e P\hbox{Diff}M$ for $e(t)=id_{M}$,
$0\le t\le T$.

 For each $0\le t\le T$ let $H_t^\tau: M\to M$, $\tau\in \R$ be the 
solution flow to the vector field $X(\cdot)(h_t)$ so
\begin{equation}\label{15}
\left\{
\begin{array}{lll}
{\partial \over\partial \tau}H_t^\tau(x)&=&X(H_t^\tau(x))h_t\\
H_t^0(x)&=&x.
\end{array}\right.
\end{equation}

\begin{lemma}
The vector field $U^h$ on $P\hbox{Diff}M$ has solution flow 
$\Phi_\tau: P \hbox{Diff}M \to P\hbox{Diff}M$, $\tau\in \R$ given by
$\Phi_\tau(\phi)_t(x)=\phi_t(H_t^\tau(x))$.
\end{lemma}

\noindent{\bf Proof.}
By left invariance we can suppose $\phi=e$. We then need only to 
observe that
$$
{\partial \over \partial \tau} H_t^\tau (x)=TH_t^\tau \left(X(x)h_t\right)$$
for each $0\le t\le T$: a standard property of ordinary,
 time-independent dynamical systems which is seen by differentiating
the identity
$$H_t^{\tau+\sigma}=H_t^\tau\circ H_t^\sigma(x)$$
with respect to $\sigma$ at $\sigma=0$.
\hfill \rule{2mm}{3mm}

\bigskip

{\bf C.} In the case where $h$ is random, with
$h: \Omega\to L_0^{2,1}([0,T];\R^d)$  adapted,  we can use the same notation
to obtain a variation of our stochastic flow $\{\xi_t: 0\le t\le T\}$
on $M$ generated by the vector field $V^h$, and given explicitly by
$$\xi^\tau_\cdot=\Phi_\tau(\xi_\cdot),$$
i.e.
\begin{equation}\label{16}
\xi_t^\tau(x)=\xi_t(H_t^\tau(x)).
\end{equation}
In particular 
\begin{equation}\label{vectorfields}
{\partial \over \partial \tau}\xi_t^\tau(x)\left |_{\tau=0}
= \right. T\xi_t\left(X(x)h_t\right).
\end{equation}

\bigskip

Using the structure of $C_x(M)$ as a $C^\infty$ Banach manifold
 let $BC^1(C_x(M))$ be the space of $C^1$ maps $F: C_x(M)\to \R$ such that
 there is a constant $|dF|_\infty$ with 
\begin{equation}
|dF(v)|\le |dF|_\infty \sup_{0\le t\le T} |v_t|
\end{equation}
for all tangent vectors $v: [0,T]\to TM $ to $C_x(M)$.
Set
$V_t^{X(h)}(x)=T\xi_t\left(X(x)(h_t)\right)$, which gives rise to a vector
field along $\{\xi_\cdot(x)\}$ on $C_x(M)$.

\begin{proposition}\label{3.2}
 Suppose $h: [0,T]\times \Omega\to T_{x}M$ is adapted, belongs
to $L_0^{2,1}$ a.s. and such that
 $\E\left(\int_0^T |\dot h_s|^2 ds\right)^{{1+\epsilon\over 2}}<\infty$
for some $\epsilon>0$.
 Then for each $x\in M$ the processes $\xi_\cdot^\tau(x)$, $\tau\in \R$ have
mutually equivalent laws ${\Bbb P}_\tau^x$, $\tau\in \R$ on
$C_x(M)$ with
$${d {\Bbb P}_\tau^x\over d {\Bbb P_0^x}}=\hbox{exp}{\left\{
\int_0^T<X(\xi_s^\tau(x))^* T\xi_s
\left({\partial \over \partial s}H_s^\tau(x)\right), dB_s>-
\half\int_0^T | T\xi_s
\left({\partial \over \partial s}H_s^\tau(x)\right)|^2 ds\right\}}.$$
Moreover, for any $F\in BC^1(C_x(M))$, 

$$\E dF( V^{X(h)}_\cdot)
=\E F(\xi_\cdot)
\int_0^T \left<  X(\xi_s(x))dB_s,
 V_s^{X(\dot h)}(x))\right>.$$
\end{proposition}

\noindent{\bf Proof.}
For the equivalent part note that $\{\xi_t^\tau: 0\le t\le T\}$
 satisfies the equation:
$$d\xi_t^\tau(x)=X\left(\xi_t^\tau(x)\right)\circ dB_t+A(\xi_t^\tau(x))dt+
T\xi_t\left({\partial \over \partial t}H_t^\tau(x)\right)dt.$$
A straightforward argument shows that 
$$\int_0^T \left|X(\xi_s^\tau(x))^* T\xi_s  \left({\partial \over \partial s}
H_s^\tau(x)\right)\right|^2\, ds <\infty, \hskip 6pt a.s.   $$
Therefore if we set
$$M_t^\tau=\int_0^t\left <X(\xi_s(x))^* T\xi_s
\left({\partial \over \partial s}H_s^\tau(x)\right), dB_s\right>,$$
then by the Girsanov-Maruyama theorem, $P_\tau^x$ is equivalent to $P_0^x$
and
\begin{equation}\label{density}
{d {\Bbb P}_\tau^x\over d {\Bbb P_0^x}}=\exp^{M_T^\tau-\half <M^\tau>_T}.
\end{equation}
Consequently,
  $$\E F(\xi_\cdot^\tau(x))=\E F(\xi_\cdot(x))
 {d {\Bbb P}_\tau^x\over d {\Bbb P_0^x}}.$$
Now suppose $h_\cdot$ and $\int_0^\cdot |\dot h_s|^2 ds$ are bounded
 on $[0,T]\times \Omega$.
Differentiating  with respect to  $\tau$  at $\tau=0$ and using
(\ref{15}) gives
$$\E dF(T\xi_\cdot (X(x)h_\cdot)) =\E F(\xi_\cdot(x))
{\partial \over \partial \tau}
\left({d {\Bbb P}_\tau^x\over d {\Bbb P_0^x}}\right)_{\tau=0},$$
since   $|dF|$ is bounded and
 $\sup_{0\le s\le T} |T\xi_s|\in \cap_{1\le  p< \infty}L^p$.

The second statement follows from differentiation of (\ref{density}),
using the fact that 
 $\left({d {\Bbb P}_\tau^x\over d {\Bbb P_0^x}}\right)_{\tau=0}=1$
and  ${\partial \over \partial t}
 H_t^\tau(x)\left|_{\tau=0}\right. =0$:

\begin{eqnarray*}
{\partial \over \partial \tau}
\left({d {\Bbb P}_\tau^x\over d {\Bbb P_0^x}}\right)_{\tau=0}
&=&\left({d {\Bbb P}_\tau^x\over d {\Bbb P_0^x}}\right)_{\tau=0}\cdot
\left[\left({\partial \over \partial \tau}M_T^\tau\right)_{\tau=0}
- \half \left({\partial \over \partial \tau}
\left <M_T^\tau \right>\right)_{\tau=0}\right] \\
&=&\int_0^T \left<X(\xi_s(x))dB_s,
 {D \over \partial \tau}\left[T\xi_s\left(
{\partial \over \partial s}H_s^\tau(x)\right)\right]\right>_{\tau=0}\\
&=&\int_0^T\left. \left<X(\xi_s(x))dB_s,
 T\xi_s({D\over \partial s} X(H_s^\tau(x))h_s)\right|_{\tau=0}
\right>\\
&=&\int_0^T \left<  X(\xi_s(x))dB_s,
 T\xi_s(X(x)\dot h_s)\right>.
\end{eqnarray*}

For general $h$ take a sequence of bounded  $h_n$ which converges to $h$ in
$L^{1+\epsilon\over 2}(\Omega, L_0^{2,1}([0,T]))$  to finish the
proof. See the proof of theorem 4.1.
\hfill\rule{2mm}{3mm}

\bigskip

The following is an analogue of Corollary \ref{2.3}: here $\tilde\nabla$
is any metric connection and $\tilde W_\cdot^Z$ is as in Corollary  \ref{2.3},

\begin{theorem}\label{3.3}
Let $F\in BC^1(C_x(M))$ and $h(\omega)\in L_0^{2,1}([0,T];\R^n)$ a.s..
Suppose
$h_\cdot$ is adapted to the filtration of $\{\F_\cdot^x\}$
and such that
 $\E\left(\int_0^T|\dot h_s|^2ds\right)^{1+\epsilon\over 2}<\infty$
 for some $\epsilon>0$.  Then

\begin{equation}
\E dF(\tilde W_\cdot^Z(h_\cdot))=     \E F(\xi_{\cdot}(x))
\int_0^T<\tilde W_s^Z(\dot h_s), \tilde{//}_sd\tilde B_s>.
\end{equation}
If  $\tilde\nabla^\prime$ is metric  for some Riemannian metric, we can
take $\epsilon=0$.
\end{theorem}

\section{Integration by parts for the free path space}

It is easy to modify the proof of Proposition \ref{3.2} to
the case where $h(0)\not =0$ and so obtain an integration
by parts formula for the free path space $PM=\cup_{x\in M} P_xM$
with uniform topology and measure given by the Riemannian measure of $M$ 
together with the laws of $\{\xi_\cdot(x): x\in M\}$.
In fact it is straightforward to generalize to the case of an
$x$-dependent $h_\cdot$. For this let $C^1(TM)$ be the space of
$C^1$ vector fields on $M$ with its usual topology:

\begin{theorem}\label{3.4}
Let $h: [0,T]\times \Omega\to   C^1(TM)$
be a cadlag adapted process such that  the $T_xM$ valued process
 $h_\cdot(x)$ has  sample paths in $L^{2,1}([0,T]; T_xM)$ for each
 $x\in M$ with  
 $|h_0(\cdot)|+ \sqrt{\int_0^t |\dot h_s(\cdot)|^2ds}$  in
$L^{1+\epsilon}\left(\Omega\times M; \R\right)$ for some $\epsilon>0$.
Let $F$ be in $ BC^1(PM;\R)$. Then
\begin{equation}\label{100}
\begin{array}{ll}
&\E \int_M  dF\left(T_x\xi_\cdot (h_\cdot (\omega)(x))\right)dx\\
=&\E \int_M  F(\xi_\cdot(x))\left\{-div h_0(x)+
\int_0^T \left<T\xi_s(\dot h_s(x)), X(\xi_s(x)) dB_s\right>\right\}dx.
\end{array}\end{equation}
\end{theorem}

\noindent{\bf Proof.}
Proceed as for Proposition \ref{3.2} but with $X(x)h_t$ replaced
by $h_t(x)$.
 In particular the definition (\ref{5}) of $H_t^\tau$ becomes

\begin{eqnarray*}
{\partial \over \partial \tau} H_t^\tau(x)&=&h_t\left(H_t^\tau(x)\right)\\
H_t^0(x)&=&x.
\end{eqnarray*}
while $\xi_t^\tau$ is defined by (\ref{16}). However now 
$\xi_0^\tau(x)=\xi_0\left(H_0^\tau(x)\right)$: the starting point is 
transported by the flow of $h_0(x)$.

 We first assume $h_\cdot$ and
$\int_0^\cdot |\dot h_s|^2ds$ are bounded on $\Omega\times M$.  
Then the   Girsanov-Maruyama  theorem gives us
equivalence between the measures $P_\tau^x$ and $P_0^{H_0^\tau(x)}$
with
$$\int_M \E F\left(\xi_\cdot^\tau(x)\right)dx
=\int_M \E F\left(\xi_\cdot(H_0^\tau(x))\right) {d{\Bbb P}_\tau^x\over
d {\Bbb P}_0^{H_0^\tau(x)}}dx.$$
On differentiating this there is the extra term
\begin{eqnarray*}
&&\int_M dF\left(T\xi_\cdot(\left.{\partial\over\partial \tau}
H_0^\tau(x)\right|_{\tau=0})\right)\, dx \\
=&&\int_M dF\left(T_x\xi_\cdot\left(h_0(x)\right)\right)dx\\
=&&\int_M d_x\left(F\circ \xi_\cdot\right)\left(h_0(x)\right)dx
\end{eqnarray*}
where $d_x\left(F\circ \xi_\cdot\right)$ refers to the derivative in $M$
of $F\circ \xi_\cdot: M\times \Omega \to \R$. Now apply the classical
 Stokes  theorem on $M$ to get:

\begin{eqnarray*}
&&\E \int_M  dF(T_x\xi_\cdot (h_\cdot(\omega)(x)))dx\\
=&&\E \int_M  F(\xi_\cdot(x))\left\{-div h_0(x)+
\int_0^T <T_x\xi_s(\dot h_s(x)), X(\xi_s(x)) dB_s>\right\}dx.
\end{eqnarray*}

 For general $h$ let $\tau_R$ be the first exit   time of
$||h_\cdot||_{C^1}+ \int_0^\cdot|h_s(x)|^2ds$  from $[0,R)$.
 Set $h_t^R(x)=h_{t\wedge \tau_R}(x)\chi_{\{||h_0||_{C^1}<R\}}$.
We have:

\begin{eqnarray*}
&& \E \int_M dF(T_x\xi_\cdot (h_\cdot^R(\omega)(x)))dx\\
=&& \E \chi_{\{||h_0||_{C^1}<R\}} \int_M   F(\xi_\cdot(x))\left\{-div h_0(x)+
\int_0^{T\wedge \tau_R} <T_x\xi_s(\dot h_s(x)), X(\xi_s(x)) dB_s>\right\}dx.
\end{eqnarray*}

Now let $R\to \infty$. The left hand side converges to
$ \E\int_M  dF(T\xi_\cdot(h_\cdot(\omega)(x)))dx$ since
$$|dF(T\xi_\cdot (h_\cdot^R(\omega)(-)))|
\le  \tilde c\sup_{t} |T\xi_t(\omega)|\sup_t|h_t(-,\omega)|$$
and 
$\sup_x \E\left(\sup_t|T\xi_t|\, \int_M \sup_{t}|h_t(x,\omega)|dx \right)<\infty$
from
\begin{eqnarray*}
&&\sup_t|h_t(x)|\le |h_0(\omega)|+\int_0^T |\dot h_s(\omega)|ds\\
&&\le |h_0(\omega)|+\sqrt T   \sqrt{\int_0^T|\dot h_s(\omega)|^2 ds}
\in L^{1+\epsilon}(\Omega\times M)
\end{eqnarray*}

 Using  Burkholder-Davis-Gundy inequality
to justify the integration on the right hand side
we see that it  converges to the right hand side of
(\ref{100}). \hfill\rule{2mm}{3mm}

\bigskip

Just as before the intrinsic formulae can be deduced using \cite{EL-LJ-LI1}:

\begin{theorem}
Let $F$ be in $BC^1(PM;\R)$ and $h$ be as in Theorem 4.1 
but with $h_\cdot(x)$ adapted  to  the filtration of $\{\F_\cdot^x\}$,
and $div h_0 \in L^1\left(\Omega\times M, \R\right)$.
 Then for any metric  connection $\tilde \nabla$ on $M$,
\begin{equation}\label{21}
\begin{array}{ll}
&\E \int_M dF\left(\tilde W^Z_\cdot (h_\cdot (\omega)(x))\right)dx\\
=&\E \int_M  F(\xi_\cdot(x))\left\{-div h_0(x)+   \int_0^T
 \left<\tilde W_s^Z(\dot h_s(x)), \tilde{//}_sd\tilde B_s\right>\right\}dx.
\end{array}\end{equation}
If furthermore $\tilde\nabla^\prime$ is metric with respect to
a Riemannian metric, we can take $\epsilon=0$.
\end{theorem}

\noindent {\bf Proof.}
The proof is just as  that of Theorem 3.3. 
\hfill\rule{2mm}{3mm}

\bigskip


\noindent
Present address of  Xue-Mei Li\\
Mathematics Department, U-9, MSB 111, University of Connecticut,\\
196 Auditorium Road, Storrs, Connecticut 06269, USA

\end{document}